\documentclass[12pt]{amsart}
\input{amssym.def}
\usepackage{graphicx}

\newtheorem{proposition}{Proposition}[section]

\theoremstyle{remark}

\newcommand{\spn}{\text{\rm span }}
\newcommand{\R}{\text{\rm Re }}

\newcommand{\Sb}{\text{\bf S}}
\newcommand{\Diff}{\text{\rm Diff }}

\newcommand{\Vect}{\text{\rm Vect }}

\newcommand{\s}{\vspace{0.3cm}}

\input epsf

\begin{document}

\title[Sub-Riemannian geometry of the coefficients...]
{Sub-Riemannian geometry of the coefficients of univalent functions}

\author[I.~Markina, D.~Prokhorov, and A.~Vasil'ev]{Irina Markina, Dmitri Prokhorov, and Alexander Vasil'ev}
\thanks{Partially supported by 
RFBR (Russia) 04-01-00083, and by the grant of the University of Bergen}
\subjclass[2000]{Primary 30C50, 17B66, 53C17; Secondary 70H06, 37J35, 81R10}
\keywords{Univalent function, coefficient, Hamiltonian system, distribution of a tangent bundle, sub-Riemannian manifold, geodesics}
\address{I.~Markina, A.~Vasil'ev: Department of Mathematics, University of Bergen, Johannes Brunsgate 12, Bergen 5008, Norway}
\email{irina.markina@uib.no}
\email{alexander.vasiliev@uib.no}
\address{D.~Prokhorov: Department of Mathematics and Mechanics, Saratov State University, Saratov 410026,
Russia}
\email{ProkhorovDV@info.sgu.ru}

\begin{abstract}  
We consider coefficient bodies $\mathcal M_n$  for univalent functions. Based on the L\"owner-Kufarev
parametric representation we get a partially integrable Hamiltonian system in which the first integrals are Kirillov's
operators for a representation of the Virasoro algebra. Then  $\mathcal M_n$  are defined as sub-Riemannian manifolds.  Given a  Lie-Poisson bracket they form a grading of subspaces
with the first subspace as a bracket-generating distribution of complex dimension two. With this sub-Riemannian structure we construct a new Hamiltonian system and calculate regular geodesics which turn to be horizontal. Lagrangian formulation is also given in the particular case $\mathcal M_3$.
\end{abstract}

\maketitle

\section{Introduction}

Let $U$ be the unit disk $U=\{z:\,\,|z|<1\}$. Let $\Sb$ stand for the standard class of holomorphic univalent functions $f:\,U \to\mathbb C$
normalized by
\[
f(z)=z\left(1+\sum\limits_{n=1}^{\infty}c_n z^n\right),\quad z\in U.
\]
By $\tilde{\Sb}$ we denote the class of functions from $\Sb$ smooth ($C^{\infty}$) on the boundary $S^1$ of $U$.
Considering $\{c_1,\dots,c_n,\dots\}$ as local affine coordinates on $\Sb$ or $\tilde{\Sb}$ we provide an embedding of these
infinite dimensional manifolds into $\mathbb C^{\infty}$. We denote by $\mathcal M$ the set $\cup_{n=1}^{\infty}\mathcal M_n$, where
\[
\mathcal M_n=\{(c_1,\dots,c_n):\,\,f\in \tilde{\Sb}\}.
\]
The class $\Sb$ is compact regarding to the local uniform topology in $U$ and $\tilde{\Sb}$($\simeq \mathcal M$) is a dense subclass of $\Sb$. By the famous de Branges' result \cite{Branges} (former Bieberbach conjecture), $\mathcal M$ lies in the bounded domain $|c_n|<n+1$, $n\geq 1$. The set $\mathcal M_1$ is the open disk $|c_1|<2$. However, the description of $\mathcal M_n$ is extremely difficult for $n>1$. Only the first non-trivial coefficient body $\mathcal M_2$ has been
described completely by Schaeffer and Spenser in 1950 in their
well-known monograph \cite{SS}. A qualitative description of $\mathcal M_n$,
$n\geq 2$, has been partially given in \cite{Babenko}. Apart from
these two monographs there are only few works where a progress in
such a problem has been made (see, e.g., \cite{Prokhorov, ProkhVas}).
Such a complicated nature of the coefficient bodies in the Euclidean structure of $\mathbb C^n$ encourages us to think of other pertinent geometries suitable to the structure of $\mathcal M_n$.

On the other hand, the manifold $\mathcal M$ is a natural representation of Kirillov's infinite dimensional K\"ahlerian manifold $\Diff S^1/S^1$ through conformal welding, here $\Diff S^1$ denotes the Lie group of orientation preserving
diffeomorphisms of the unit circle $S^1$, and the subgroup of rotations is associeated with $S^1$. Indeed, given a map $f\in \tilde\Sb$ we construct an adjoint univalent
meromorphic map 
\[
g(z)=d_1z+d_0+\frac{d_{-1}}{z}+\dots,
\] 
defined in the exterior $U^*$ of $U$, and such that  $\hat{\mathbb C}\setminus\overline{f(U)}=g(U^*)$. This gives the identification $\Diff S^1/S^1$ with $\mathcal M$, see \cite{Airault, KY1}. The central extension of $\Diff S^1$ by $\mathbb R$ is the Virasoro-Bott group. The corresponding central extension of the space $\Vect S^1$ of vector fields on $S^1$ is the Virasoro algebra ($= \Vect S^1\oplus \mathbb R$). The infinitesimal action of $\Diff S^1$ on $\mathcal M$ (given by the Goluzin-Schiffer variation) leads to special vector fields $L_j$ on $\mathcal M$, Kirillov's
operators for a representation of the Virasoro algebra. 

We deduce a Hamiltonian system for the L\"owner-Kufarev trajectories in $\mathcal M_n$. In view of Hamiltonian mechanics, this formulation performs a trivial motion with constant speed (and vanishing energy).
Our aim is to describe a sub-Riemannian structure of the $n$-complex-dimensi\-onal manifold $\mathcal M_n$ based on Kirillov's operators and to describe geodesics in this structure. We calculate them explicitly for $n=3$.  Such a description gives a non-trivial motion in which the energy of the system conserves along non-Riemannian geodesics.

In our setup Kirillov's operators appear as the first integrals
of a partially integrable Hamiltonian system for $c_n$ generated by the L\"owner-Kufarev representation of univalent functions. The sub-Riemannian structure is based on the distribution defined by only two first vector fields $L_1$ and $L_2$ and other vector fields form a grading sequence. The horizontal curves are only of finite length in the corresponding sub-Riemannian metric and we give a  description of regular geodesics in $\mathcal M$. Lagrangian formulation is also given in the particular case $\mathcal M_3$.

\section{Hamiltonian system for the coefficients}

\subsection{Coefficient bodies} By the {\it coefficient problem for univalent functions} we mean the
problem of precise finding the regions $\mathcal M_n$ defined above. These
sets have been investigated by a great number of authors,
but the most remarkable source is a monograph \cite{SS}
written by Schaeffer and Spencer in 1950. Among other contributions
to the coefficient problem we distinct a monograph by Babenko
\cite{Babenko} that contains a good collection of qualitative
results on the coefficient bodies $\mathcal M_n$. The results concerning the
structure and properties of $\mathcal M_n$ include (see \cite{Babenko}, \cite{SS})

\begin{itemize}

\item[(i)] $\mathcal M_n$ is homeomorphic to a $(2n-2)$-dimensional ball and its
boundary $\partial \mathcal M_n$ is homeomorphic to a $(2n-3)$-dimensional
sphere;

\item[(ii)] every point $x\in \partial \mathcal M_n$ corresponds to exactly one function $f\in
\Sb$ which is called a {\it boundary function} for $\mathcal M_n$;

\item[(iii)] with the exception for a set of smaller dimension,
at every point $x\in \partial \mathcal M_n$ there exists a normal vector
satisfying the Lipschitz condition;

\item[(iv)] there exists a connected open set $X_1$ on $\partial \mathcal M_n$,
such that the boundary $\partial \mathcal M_n$ is an analytic hypersurface at
every point of $X_1$. The points of $\partial \mathcal M_n$ corresponding to
the functions that give the extremum to a linear functional belong
to the closure of $X_1$.

\end{itemize}

It is worth to note again that all boundary functions have a similar
structure. They map the unit disk $U$ onto the complex plane
$\mathbb C$ minus piecewise analytic Jordan arcs forming a tree with
a root at infinity and having at most $n$ tips, as it has been mentioned in the
preceding section.  The uniqueness of the boundary functions
implies that each point  of $\partial \mathcal M_n$ (the set of first coefficients) defines
the rest of coefficients uniquely.  

\subsection{Hamiltonian dynamics and integrability}

Let us recall briefly the Hamiltonian and symplectic definitions and
concepts that will be used in the sequel. There exists a vast amount
of modern literature dedicated to different approaches to and
definitions of {\it integrable systems} (see, e.g., \cite{Arnold},
\cite{Babelon}, \cite{Bolsinov}, \cite{Zakharov}).

The classical definition of a {\it completely integrable system} in
the sense of Liouville applies to a Hamiltonian system. If we can
find  independent conserved integrals which are pairwise involutory
(vanishing Poisson bracket), this system is completely integrable
(see e.g., \cite{Arnold}, \cite{Babelon}, \cite{Bolsinov}). That is
each first integral allows us to reduce the order of the system not
just by one, but by two. We formulate this definition in a slightly
adopted form as follows.

A dynamical system in $\mathbb C^{2n}$ is called {\it Hamiltonian}
if it is of the form
\begin{equation}
\dot{x}=\nabla_s H(x),\label{Hamilton1}
\end{equation}
where $\nabla_s$ denotes the {\it symplectic gradient} given by
\[
\nabla_s=\left(\frac{\partial}{\partial
\bar{x}_{n+1}},\dots,\frac{\partial}{\partial
\bar{x}_{2n}},-\frac{\partial}{\partial
{x}_1},\dots,-\frac{\partial}{\partial {x}_n}\right).
\]
The function $H$ in (\ref{Hamilton1}) is called the {\it Hamiltonian
function} of the system. It is convenient to redefine the
coordinates $(x_{n+1},\dots, x_{2n})=(\psi_{1},\dots, \psi_{n})$,
and rewrite the system (\ref{Hamilton1}) as
\begin{equation}
\dot{x}_k=\frac{\partial H}{\partial \overline{\psi}_k},\quad
\dot{\overline{\psi}}_k=-\frac{\partial H}{\partial x_k},\quad
k=1,2\dots,n.\label{Hamilton2}
\end{equation}
The system has $n$ degrees of freedom. The two-form
$\omega=\sum_{k=1}^n dx\wedge d\bar{\psi}$ admits the 
Lie-Poisson bracket $[\cdot,\cdot]$
\[
[f,g]=\sum\limits_{k=1}^{n}\left(\frac{\partial f}{\partial x_k}
\frac{\partial g}{\partial \overline{\psi}_k}-\frac{\partial
f}{\partial \overline{\psi}_k} \frac{\partial g}{\partial
x_k}\right)
\]
associated with $\omega$. The symplectic pair $(\mathbb C^{2n},
\omega)$ defines the Poisson manifold $(\mathbb C^{2n},
[\cdot,\cdot])$. These notations may be generalized for a
symplectic manifold and a Hamiltonian dynamical system on it.

The system (\ref{Hamilton2}) may be rewritten  as
\begin{equation}
\dot{x}_k=[x_k, H],\quad
\dot{\overline{\psi}}_k=[\overline{\psi}_k, H],\quad
k=1,2\dots,n,\label{Hamilton3}
\end{equation}
and the {\it first integrals} $L$ of the system are characterized
by
\begin{equation}
[L, H]=0.\label{Hamilton4}
\end{equation}
In particular, $[H,H]=0$, and the Hamiltonian function $H$ is an
integral of the system (\ref{Hamilton1}). If the system
(\ref{Hamilton3}) has $n$ functionally independent integrals
$L_1,\dots,L_n$, which are pairwise involutory
$[L_k,L_j]=0$, $k,j=1,\dots,n$, then it is called {\it
completely integrable} in the sense of Liouville. The function $H$
is included in the set of the first integrals. The classical theorem
of Liouville and Arnold \cite{Arnold} gives a complete description
of the motion generated by the completely integrable system
(\ref{Hamilton3}). It states that such a system admits action-angle
coordinates around a connected regular compact invariant manifold.

If the Hamiltonian system admits only $1\leq k<n$ independent
involutory integrals, then it is called {\it partially integrable}.
The case $k=1$ is known as the Poincar\'e--Lyapunov theorem which
states that a periodic orbit of an autonomous Hamiltonian system can
be included in a one-parameter family of such orbits under a
non-degeneracy assumption. A bridge between these two extremal cases
$k=1$ and $k=n$ has been proposed by Nekhoroshev \cite{Nekhoroshev}
and proved later in \cite{Bambusi}, \cite{Fiorani}, \cite{Gaeta}.
The result states the existence of $k$-parameter families of tori
under suitable non-degeneracy conditions.

\subsection{Hamiltonian system for the coefficients}
The L\"owner-Kufarev parametric method (see, e.g., \cite{Duren, Pom}) is based on a representation of any function $f$ from the class $\Sb$ by the limit
\begin{equation}
f(z)=\lim\limits_{t\to\infty}e^tw(z,t)\label{lim},
\end{equation}
where the function $$w(z,t)=e^{-t}z\left(1+\sum\limits_{n=1}^{\infty}c_n(t) z^n\right)$$
is a solution to the L\"owner-Kufarev equation
\begin{equation}
\frac{dw}{dt}=-wp(w,t),\label{LK}
\end{equation}
with the initial condition $w(z,0)\equiv z$. The function $p(z,t)=1+p_1(t)z+\dots$ is holomorphic in $U$ and has the positive real part for all $z\in U$ almost everywhere in $t\in [0,\infty)$. If $f\in \tilde\Sb$, then
\begin{eqnarray}
\dot{c}_n & = &c_n-\frac{e^t}{2\pi i}\int\limits_{S^1}w(z,t)p(w(z,t),t)\frac{dz}{z^{n+2}},\label{coeff}\\ 
&=&-\frac{1}{2\pi i}\int\limits_{S^1}\sum\limits_{k=1}^ne^{-kt}(e^tw)^{k+1}p_k\frac{dz}{z^{n+2}}\quad n\geq 1.\nonumber
\end{eqnarray}
In particular,
\begin{eqnarray*}
\dot{c}_1 & = & -e^{-t}p_1,\\
\dot{c}_2 & = & -2e^{-t}p_1c_1-e^{-2t}p_2,\\
\dot{c}_3 & = & -e^{-t}p_1(2c_2+c_1^2)-3e^{-2t}p_2c_1-e^{-3t}p_3,\\
\dots& & \dots
\end{eqnarray*}

We consider an adjoint vector
\[
\psi(t)=\left(
\begin{array}{c}
\psi_1(t)\\
\cdot\\
\cdot\\
\cdot\\
\psi_n(t)
\end{array} \right),
\]
with complex-valued coordinates $\psi_1,\dots,\psi_n$, and the
complex Hamiltonian function
\[
H(a,{\psi},u)=\sum\limits_{k=1}^n\bar{\psi}_k\left(c_k-\frac{e^t}{2\pi i}\int\limits_{S^1}w(z,t)p(w(z,t),t)\frac{dz}{z^{k+2}}\right).
\]
 To come to the Hamiltonian formulation for the
coefficient system we require that $\bar{\psi}$ satisfies the
adjoint to (\ref{coeff}) system of differential equations
\begin{equation*}
\dot{\bar{\psi}}_j=-\frac{\partial H}{\partial c_j},\quad 0\leq
t<\infty,
\end{equation*}
or
\begin{equation}
\dot{\bar{\psi}}_j=-\bar{\psi}_j+\frac{1}{2\pi i}\sum\limits_{k=1}^n\bar{\psi}_k\int\limits_{S^1}(p+wp')\frac{dz}{z^{k-j+1}},\quad j= 1,\dots, n-1,\label{psi1}
\end{equation}
and 
\begin{equation}
\dot{\bar{\psi}}_n=0.\label{psi2}
\end{equation}
In particular, for $n=3$ we have
\begin{eqnarray*}
\dot{\bar{\psi}}_1 & = & 2e^{-t}p_1\bar{\psi}_2+(2e^{-t}p_1c_1+3e^{-2t}p_2)\bar{\psi}_3,\\
\dot{\bar{\psi}}_2 & = & 2e^{-t}p_1\bar{\psi}_3,\\
\dot{\bar{\psi}}_3 & = & 0.
\end{eqnarray*}

\subsection{First integrals and partial integrability}
Let us construct the following series
\begin{equation}
\sum\limits_{k=1}^n\bar{v}_{n-k+1} z^{k-1}=e^tw'(z,t)\sum\limits_{k=1}^n\bar{\psi}_{n-k+1} z^{k-1}+e^tw'(z,t)\sum\limits_{k=n}^{\infty}b_k z^k.\label{firstint}
\end{equation}
Taking into account (\ref{psi1}) and the formula for the derivative 
\[
\frac{\partial (e^tw')}{\partial t}=e^tw'(1-p(w,t)-wp'(w,t)),
\]
we come to the conclusion that $\dot{\bar v}=0$ and $\bar v$ is constant.
We denote by  $(L_1,\dots, L_n)^T$ the vector of the
first integrals of the Hamiltonian system (\ref{coeff}-- \ref{psi2})
given by
\begin{equation}
\left(
\begin{array}{c}
L_1\\
L_2\\
L_3\\
\dots\\
L_n
\end{array} \right)=\left(
\begin{array}{rrrrr}
1&2c_1 & \dots & (n-1)c_{n-2} & nc_{n-1}\\
0 & 1 & \dots &  (n-2)c_{n-3} & (n-1)c_{n-2}\\
0 & 0 & \dots & (n-3)c_{n-4} & (n-2)c_{n-3}\\
\dots & \dots & \dots & \dots & \dots\\
0 & 0 & \dots & 0 & 1
\end{array} \right)\left(
\begin{array}{c}
\bar{\psi}_1\\
\bar{\psi}_2\\
\bar{\psi}_3\\
\dots\\
\bar{\psi}_n
\end{array} \right).\label{L1}
\end{equation}
Indeed, the
equality (\ref{firstint}) implies that $L_k=\bar{v}_k$ are
constants for all $t$ and $k=1,\dots,n$. Naturally, 
\[
[L_j, H]=\sum\limits_{k=1}^n\frac{\partial L_j}{\partial c_k}\frac{\partial H}{\partial \overline{\psi_k}}-\frac{\partial L_j}{\partial \overline{\psi_k}}\frac{\partial H}{\partial c_k}= \sum\limits_{k=1}^n\frac{\partial L_j}{\partial c_k}\dot{c}_k+\frac{\partial L_j}{\partial \overline{\psi_k}}\dot{\bar{\psi}}_k=\dot{L}_j=0.
\]
The commutator relations are:
\begin{equation}
[L_j, L_k]=(j-k)L_{k+j}, \quad \mbox{when $k+j\leq n$,}\label{L3}
\end{equation}
or 0 otherwise. This implies that
\begin{itemize}
\item  the first integrals $(L_{[(n+1)/2]},\dots, L_n)$ are
pairwise involutory;
\item the integrals $(L_1,\dots,L_{[(n-1)/2]})$ are not pairwise
involutory but their Lie-Poisson brackets give all the rest of
integrals. 
\end{itemize}
It is clear from the form of the matrix in the above representation of $L_k$, $k=1,\dots, n,$
that all these integrals are algebraically (even linearly) independent. Therefore, the Hamiltonian system (\ref{coeff}--\ref{psi2}) is partially integrable in the Liouville sense. In particular for $n=3$, we compute 
\begin{eqnarray*}
\psi_1&=&(4c_1^2-3c_2)v_3-2c_1v_2+v_1,\\
\psi_2&=&-2c_1v_3+v_2,\\
\psi_3&=&v_3.
\end{eqnarray*}
\s

\noindent{\it Remark.} All previous considerations we did for the class $\tilde{\Sb}$ because it will be important for us
in the future sections. But the result on partial integrability is still valid for the whole class $\Sb$ going inside
the unit disk by $f\to \frac{1}{r}f(rz)$, and letting $r\to 1$.

\s

\noindent{\it Remark.} The complete integration of this Hamiltonian system requires additional information
on the trajectories, in particular, on the controls $p_1,p_2,\dots$. One way to perform such integration is
solution of the extremal problem of finding the boundary hypersurfaces of $\mathcal M_n$ by optimal control methods, see \cite{ProkhVas}.

\s

\noindent{\it Remark.} In view of Hamiltonian mechanics, our Hamiltonian system describes  trivial motion with the constant velocity because the Hamiltonian function is linear with respect to $\psi$. An attempt to get a non-trivial
description of the L\"owner-Kufarev motion was launched in \cite{Vas} by intaking a special Lagrangian. Further on in this paper, we shall
give another non-trivial Hamiltonian and Lagrangian descriptions based on the sub-Riemannian geometry led on $\mathcal M_n$.

\s

\noindent{\it Remark.} The coefficient bodies $\mathcal M_1$, $\mathcal M_2, \dots$ generate a hierarchy of Hamiltonian systems 
(\ref{coeff}-\ref{psi1}).

\section{Virasoro algebra and Kirillov's operators}

A {\it Killing vector field} is a vector field on a Riemannian manifold that preserves the metric. Killing fields are the infinitesimal generators of isometries; that is, flows generated by Killing fields are continuous isometries of the manifold. A {\it Witt algebra} is the Lie algebra of Killing vector fields defined on the Riemann sphere. The basis for these Killing fields is given by the holomorphic fields 
\[
L_n=-z^{n+1} \frac{\partial}{\partial z}. 
\]
The Lie-Poisson bracket of two Killing fields is
\begin{equation}
[L_m,L_n]=(n-m)z^{m+n+1}\frac{\partial}{\partial z}=(m-n)L_{m+n}.\label{Lie1}
\end{equation}
The {\it Virasoro algebra} is the central extension of the Witt algebra by $\mathbb C$. The Lie-Poisson bracket for the basis vectors of the Virasoro algebra is
\[
[L_m,L_n]_{Vir}=(m-n)L_{m+n}+\frac{c}{12}n(n^2-1)\delta_{n,-m}.
\]
The constant $c\in \mathbb C$ is known as the {\it central charge} and is a constant of the theory.

To analyze and to represent this central extension we 
consider real vector fields over the unit circle. We denote the Lie group of $C^{\infty}$ sense preserving
diffeomorphisms of the unit circle $S^1$  by $\Diff S^1$. Each element of $\Diff S^1$ is represented as
$z=e^{i\alpha(\theta)}$ with a monotone increasing $C^{\infty}$
real-valued function $\alpha(\theta)$, such that
$\alpha(\theta+2\pi)=\alpha(\theta)+2\pi$.
 The Lie algebra for $\Diff S^1$ is
identified with the Lie algebra $\Vect S^1$ of smooth ($C^{\infty}$)
tangent vector fields to $S^1$, the infinitesimal action is $\theta\to\theta+\varepsilon \phi(\theta)$. To $\phi$ we associate the vector field $\phi\frac{d}{d\theta}$, and the Lie-Poisson 
bracket is given by
$$[\phi_1,\phi_2]={\phi}_1{\phi}'_2-{\phi}_2{\phi}'_1. $$
Fixing  the trigonometric basis in $\Vect S^1$, the commutator
relations admit the form
\begin{eqnarray*}
\left[\cos\,n\theta, \cos\,m\theta\right] & = &
\frac{n-m}{2}\sin\,(n+m)\theta+ \frac{n+m}{2}\sin\,(n-m)\theta,\\
\left[\sin\,n\theta, \sin\,m\theta\right]& = &
\frac{m-n}{2}\sin\,(n+m)\theta+ \frac{n+m}{2}\sin\,(n-m)\theta,\\
\left[\sin\,n\theta, \cos\,m\theta\right] & = &
\frac{m-n}{2}\cos\,(n+m)\theta- \frac{n+m}{2}\cos\,(n-m)\theta.
\end{eqnarray*}
The space $\Vect S^1$ with so given Lie bracket is the space of left-invariant vector fields.

Let $I$ and $G$ be Lie algebras.  An exact sequence is a sequence of objects and morphisms between them, such that the image of one morphism equals the kernel of the next. Let us consider the exact sequence of Lie algebras 
\[
0
\longrightarrow I\stackrel{f}{\longrightarrow}
E\stackrel{g}{\longrightarrow}
G\longrightarrow 0.
\]
$E$ is called the {\it central extension} of $G$ by $I$ if $I$ belongs to the center of $E$. The central extension is given as $E\simeq G\oplus I$. A simple example is $[x+a]_{E}=[x,y]_{G}+[a,b]_{I}$.
The (real) Virasoro algebra is the unique (up to isomorphism) non-trivial central extension of $\Vect\, S^1$ by $\mathbb R$ given by the {\it Gelfand-Fuchs cocycle}~\cite{Gelfand}:
\[
\omega(\phi_1,\phi_2)=\frac{1}{2\pi}\int\limits_{0}^{2\pi}(\phi'_1\phi''_2-\phi''_1\phi'_2)d\theta.
\]
The Virasoro algebra $Vir$ is a Lie algebra over the space $\Vect\,S^1\oplus \mathbb R$ defined by the commutator
\[
[(\phi_1,a),(\phi_2, b)]_{Vir}=([\phi_1,\phi_2]_{\Vect\,S^1},\frac{c}{12}\omega(\phi_1,\phi_2)),
\]
where $a$ and $b$ are elements of the center, $ab-ba$ vanishes, and $c\in \mathbb R$ is the central charge. Integration by parts leads to
the 2-cocycle condition
\[
\omega(\phi_1,[\phi_2,\phi_3])+\omega(\phi_2,[\phi_3,\phi_1])+\omega(\phi_3,[\phi_1,\phi_2])=0,
\]
and
\begin{equation}
\omega(\phi_1,\phi_2)=-\frac{1}{4\pi}\int\limits_{0}^{2\pi}(\phi'_1+\phi'''_1)\phi_2d\theta.\label{cocycle}
\end{equation}
Correspondingly, we consider the group $\Diff\, S^1$. The {\it Virasoro-Bott group} is the unique (up to isomorphism) non-trivial central extension of $\Diff\, S^1$ given by the Thurston-Bott cocycle \cite{Bott}
\[
\Omega(f,g)=\frac{1}{2\pi}\int\limits_{0}^{2\pi}\log((f\circ g)')d\log (g').
\]
The Virasoro-Bott group is given by the following product on $\Diff\,S^1\times \mathbb R$
\[
(f,\alpha)(g,\beta)=(f\circ g, \alpha+\beta+\frac{c}{12}\Omega(f,g)).
\]

We shall identify $\Vect\, S^1$  with the functions with vanishing mean value over~$S^1$. It gives
\[
\phi(\theta)=\sum\limits_{n=1}^{\infty}a_n\cos\,n\theta+b_n\sin\,n\theta.
\]
Let us define a complex structure by the operator
\[
J(\phi)(\theta)=\sum\limits_{n=1}^{\infty}-a_n\sin\,n\theta+b_n\cos\,n\theta.
\]
Then $J^2=-id$. On $\Vect\,S^1\oplus \mathbb C$, the operator $J$ diagonalizes and we have
\[
\phi\to \phi-iJ(\phi)=\sum\limits_{n=1}^{\infty}(a_n-ib_n)e^{in\theta},
\]
and the latter extends into the unit disk as a holomorphic function.

Taking the basis of $\Vect\,S^1\oplus \mathbb C$ as $e_n=-ie^{in\theta}\partial$ we get
\[
[e_n,e_m]=(n-m)e_{n+m}+\frac{c}{12}n(n^2-1)\delta_{n,-m}.
\]
The Virasoro algebra is realizable both as a central extension of the Witt algebra and as an algebra of the Virasoro generators in Conformal Field Theory.

There is no general theory of infinite dimensional Lie groups,
example of which is under consideration.  The interest to the
particular case $\Diff S^1$ comes first of all from the two-dimensional Conformal Field Theory where the
algebra of energy momentum tensor deformed by a central extension due to the conformal anomaly
is represented by the Virasoro algebra. Entire necessary background for the
construction of the theory of unitary representations of $\Diff S^1$
is found in the study of Kirillov's homogeneous K\"ahlerian manifold
$\Diff S^1/S^1$. The group $\Diff S^1$ acts as a group of
translations on the manifold $\Diff S^1/S^1$ with $S^1$ as a stabilizer.
The K\"ahlerian geometry of $\Diff S^1/S^1$ has been described by Kirillov and
Yuriev in \cite{KY1}. The manifold $\Diff S^1/S^1$ admits several
representations, in particular, in the space of smooth probability
measures, symplectic realization in the space of quadratic
differentials. We shall use its analytic representation by $\tilde \Sb$ and $\mathcal M$ mentioned in Introduction.

The Kirillov infinitesimal action of $\Vect S^1$ on $\tilde \Sb$ is given by
the Goluzin-Schiffer variational formulas which lift the actions from the Lie algebra $\Vect
S^1$ onto $\tilde\Sb$. Let $f\in\tilde{\Sb}$ and let
$\nu(e^{i\theta})$ be a $C^{\infty}$ real-valued function in
$\theta\in(0,2\pi]$ from $\Vect S^1$ making an infinitesimal action
as $\theta \mapsto \theta+\varepsilon \nu(e^{i\theta})$. Let us consider a
variation of $f$ given by
\begin{equation}
\delta_{\nu}f(z)=\frac{f^2(z)}{2\pi
i}\int\limits_{S^1}\left(\frac{wf'(w)}{f(w)}\right)^2\frac{\nu(w)dw}{w(f(w)-f(z))} .\label{var}
\end{equation}
 Kirillov and Yuriev \cite{KY1}, \cite{KY2} (see also \cite{Airault}) have established
that the variations $\delta_{\nu}f(\zeta)$ are closed with respect to
the commutator (\ref{Lie1}) and the induced Lie algebra is the same as $\Vect
S^1$. Moreover, Kirillov's result \cite{Kir} states that there is
an exponential map $\Vect S^1\to \Diff S^1$ such that the subgroup
$S^1$ coincides with the stabilizer of the map $f(z)\equiv
z$ from $\tilde{\Sb}$.

Taking the complexification $\Vect_{\mathbb C} S^1$ of $\Vect S^1$ and the basis $\nu=-iz^k$ in the integrand of (\ref{var}) we calculate the residue in (\ref{var}) and obtain
\[
L_k(f)(z)=\delta_{\nu}f(z)=z^{k+1}f'(z), \quad k=1,2,\dots
\]
In terms of the affine coordinates in $\mathcal M$ we get
\[
L_j=\partial_j+\sum\limits_{k=1}^{\infty}(k+1)c_{k}\partial_{j+k},
\]
or truncating
\begin{equation}
L_j=\partial_j+\sum\limits_{k=1}^{n-j}(k+1)c_{k}\partial_{j+k},\label{L2}
\end{equation}
on $\mathcal M_n$, where $\partial_k=\partial/\partial c_k$. Considering the adjoint vector $\psi$ (Section 2) as the vector
of affine coordinates, we  conclude that the vector fields given by the first integrals $L_k$, see (\ref{L1}), are exactly  Kirillov's operators. Given a fixed central charge $c$, Neretin \cite{Neretin} introduced the sequence of polynomials $P_n$ defined by the following recurrence relations
\[
L_k(P_j)=(j+k)P_{j-k}+\frac{c}{12}k(k^2-1)\delta_{j,k},\quad P_0\equiv P_1\equiv 0,\,\,\,P_j(0)=0.
\]
Representing the momentum-energy tensor in the 2-D Conformal Field Theory the Schwarzian derivative naturally comes into play in the definition of $P_n$. It turns out that
\[
\frac{cz^2}{12}S_f(z)=\sum\limits_{n=0}^{\infty}P_n(c_1,\dots, c_n)z^n,
\]
where
\[
S_f(z)=\frac{f'''(z)}{f'(z)}-\frac{3}{2}\left(\frac{f''(z)}{f'(z)}\right)
\]
is the Schwarzian derivative of a univalent function $f\in \tilde\Sb$. In particular,
\[
\frac{1}{c}P_2(c_1,c_2)=\frac{1}{2}(c_2-c_1^2),\quad \frac{1}{c}P_3(c_1,c_2,c_3)=2(c_3-2c_1c_2+c_1^3),
\]
\[
\frac{1}{c}P_4(c_1,c_2,c_3,c_4)=5c_4-10c_1c_3-6c_2^2+17c_1^2c_2-6c_1^4,\dots
\]
\s

\noindent{\it Remark.} In general, we have real vector fields in $\Vect S^1$. The computation of $L_k$ must be carried out with respect to the basis $1,e^{\pm k i\theta}$ that leads also to $L_k$ with $k\leq 0$. However, we deal with holomorphic functions and $L_k$ with $k> 0$ are to be treated as complex vector fields (see discussion in \cite[p. 738]{Kir2}, \cite[p. 632--634]{Airault}).

\section{Sub-Riemannian geometry of $\mathcal M_n$}

A {\it sub-Riemannian structure} on an $n$-dimensional manifold $\mathcal M_n$ is a smoothly varying distribution $\mathcal D$ of $k$-planes together with a smoothly varying scalar product on these planes. The distribution $\mathcal D$ is a linear sub-bundle of a tangent bundle $T\mathcal M_n$ of $\mathcal M_n$. The {\it dimension} of the sub-Riemannian manifold is the pair $(k,n)$ (see, e.g., \cite{Montgomery, Strichartz1, Strichartz2}). In the case $n=k$ we come to the standard Riemannian structure. If $k<n$, then several new phenomena occur, e.g., the Hausdorff dimension is larger than the topological dimension,  the space of paths joining two fixed points and tangent to the distribution can have singularities. Suppose that a system of vector fields $X_1,\dots, X_k$ form an orthonormal basis of $\mathcal D$ with respect to an inner product $\langle\cdot,\cdot\rangle$. The pair $(\mathcal D,\langle\cdot,\cdot\rangle)$ is called a {\it sub-Riemannian metric} on $\mathcal M_n$. A {\it horizontal path} is an absolutely continuous path $\gamma:\,[0,1]\to \mathcal M_n$ with a tangent vector $\dot{\gamma}$
in $\mathcal D$: i.e., $\dot{\gamma}(t)=\sum_{j=1}^ku_j(t)X_j(\gamma(t))$. The length of such a path is 
\[
\int_{[0,1]}\sqrt{\langle\dot{\gamma}(t),\dot{\gamma}(t)\rangle}dt.
\]
The distance between two points is the infimum of the length of horizontal curves joining them. It is called the {\it Carnot-Carath\'eodory distance} in the literature (e.g., \cite{NSW}). Sub-Riemannian structures appear in the works of Carnot on thermodynamics and Carath\'eodory was inspired by his ideas. If all vector fields $X_1,\dots, X_k$ together with
their commutators form the total tangent space $T\mathcal M_n$, then is said that $X_1,\dots, X_k$ satisfy the bracket generating condition (or H\"ormander's hypoellipticity condition \cite{Hoermander}). The number of thee commutators is independent of the point of $\mathcal M_n$. If the manifold $\mathcal M_n$ is connected (what is satisfied in our case), and the bracket generating condition holds, then any two points can be connected by a smooth horizontal path \cite{Chow, Rashevski}.

\subsection{Sub-Riemannian structure defined by Kirillov's operators}

\begin{proposition}
Let $\mathcal M_n$ be the $n$-th coefficient body and $L_1,\dots L_n$ be vector fields defined by (\ref{L1}, \ref{L2}).
Then the system $(L_1,L_2)$ satisfies the bracket generating condition and the distribution is $\mathcal D=\spn(L_1,L_2)$.
\end{proposition}

\begin{proof}
The commutator relations (\ref{L3}) imply that the vector field $L_3$ is a unique vector generated by $L_1$ and $L_2$ by $[L_2,L_1]=L_3$. We denote by $\mathcal D_1$ the vector space generated by $L_3$. By $\mathcal D_k$ we denote the vector space given by the recurrence process $\mathcal D_k=[\mathcal D,\mathcal D_{k-1}]\setminus \mathcal D_{k-1}$. Thus, $\mathcal D_2=\spn(L_4,L_5)$, $\mathcal D_3=\spn(L_6,L_7)$. For even $n$ we have the last space $\mathcal D_{n/2}=\spn (L_n)$. For odd $n$ the last space is $\mathcal D_{(n-1)/2}=\spn (L_{n-1},L_n)$. The vector spaces 
\[
\mathcal D\oplus \mathcal D_1\oplus\dots \oplus\mathcal D_{[n/2]}=T\mathcal M_n,
\] 
form a grading sequence in $T\mathcal M_n$. The number $[n/2]$ is the degree of non-holonomy. Obviously, given $L_1,L_2$ we construct all other vector fields $L_k$, $k=3,\dots, n,$  by recurrence of commutators and 
\[
T\mathcal M_n=\spn(L_1,\dots, L_n).
\]
\end{proof}

The scalar product $\langle \cdot,\cdot\rangle$ on $\mathcal D$ will be defined by the K\"ahlerian structure of $\mathcal M_n$.  Thus, the triple $(\mathcal M_n, \mathcal D, \langle \cdot,\cdot\rangle)$ is a sub-Riemannian manifold. By abuse of notation, let us denote it simply by $\mathcal M_n$.

\begin{proposition}
The Hausdorff (complex) dimension of the sub-Riemannian manifold $\mathcal M_n$ is equal to
\begin{itemize}
\item $(\frac{n}{2}+1)^2-\frac{9}{4}$ for odd $n$;
\item $(\frac{n}{2}+1)^2-2$ for even $n$.
\end{itemize}
\end{proposition}
\begin{proof}
Let us consider the case of odd $n$. 
The complex topological dimension $\dim_{\mathbb C}\mathcal D=2$, $\dim_{\mathbb C}\mathcal D_1=1$, $\dim_{\mathbb C}\mathcal D_k=2$, for $k\geq 2$. The following formula \cite{Mitchell, Pansu} is used to calculate the Hausdorff dimension
of $\mathcal M_n$:
\[
\dim_{\mathbb C}\mathcal D+2\dim_{\mathbb C}\mathcal D_1+3\dim_{\mathbb C}\mathcal D_2+\dots+(\frac{n-1}{2}+1)\dim_{\mathbb C}\mathcal D_{\frac{n-1}{2}}=(\frac{n}{2}+1)^2-\frac{9}{4}.
\]
For even $n$ we observe that the dimension of the last subspace is 1.
\end{proof}

\begin{proposition}\label{hcprop}
A path $\gamma(s)=(c_1(s),\dots,c_n(s))$ in $\mathcal M_n$ is horizontal if and only if
\begin{equation}
\label{hc}
\begin{split}
& \dot c_3(s) =3c_2(s)\dot c_1(s)+2c_1(s)\big(\dot c_2(s)-2c_1(s)\dot c_1(s)\big)\\ & \ldots\ldots\\ & \dot c_n(s)=nc_{n-1}(s)\dot c_1(s)+(n-1)c_{n-2}(s)\big(\dot c_2(s)-2c_1(s)\dot c_1(s)\big).
\end{split}
\end{equation}
\end{proposition}

\begin{proof}
The tangent vector to $\gamma(s)$ in the local affine basis $\partial_1,\ldots,\partial_n$ is  $$\dot\gamma(s)=\dot c_1(s)\partial_1+\ldots+\dot c_n(s)\partial_n.$$ Let us rewrite the  tangent vector $\dot\gamma(s)$ in the local basis $L_1,L_2$ of the distribution $\mathcal D$. We get 
\begin{eqnarray*}
\dot\gamma(s) & =&\dot c_1(s)\partial_1+\ldots+\dot c_n(s)\partial_n\\ 
&= &\dot c_1(s)(\partial_1+2c_1\partial_2+\ldots+nc_{n-1}\partial_n)\\ 
& &+(\dot c_2(s)-2c_1\dot c_1)(\partial_2+2 c_1\partial_3+\ldots+(n-1)c_{n-2}\partial_n)\\ 
& &- \dot c_1(s)(2c_1\partial_2+\ldots+nc_{n-1}\partial_n)\\ 
& &-(\dot c_2(s)-2c_1\dot c_1)(2c_1\partial_3+\ldots+(n-1)c_{n-2}\partial_n)\\ 
& &+2c_1\dot c_1\partial_2+\dot c_3(s)\partial_3+\ldots+\dot c_n(s)\partial_n\\ 
& =&\dot c_1(s)L_1(\gamma(s))+(\dot c_2(s)-2c_1(s)\dot c_1(s))L_2(\gamma(s))\\
& &+(\dot c_3(s)-3c_2(s)\dot c_1(s)-2 c_1(s)(\dot c_2-2c_1\dot c_1))\partial_3+\ldots\\
& & +(\dot c_n(s)-nc_{n-1}\dot c_1(s)-(n-1)c_{n-2}(\dot c_2(s)-2c_1\dot c_1))\partial_n.
\end{eqnarray*} 

To simplify the calculations we use the following notation $u_1=\dot c_1$, $u_2=\dot c_2(s)-2c_1(s)\dot c_1(s)$, and $g_{k}=\dot c_{k}(s)-kc_{k-1}(s)\dot c_1(s)-(k-1)c_{k-2}(s)(\dot c_2-2c_1\dot c_1)=\dot c_{k}(s)-kc_{k-1}(s)u_1-(k-1)c_{k-2}(s)u_2$. Then  \begin{equation*}
\begin{split}
\dot\gamma(s) & =u_1L_1+u_2L_2+g_3L_3+(-2g_3c_1+g_4)\partial_4+\ldots+(-(n-2)g_3c_{n-3}+g_n)\partial_n.
\end{split}\end{equation*} Since the path $\gamma$ is supposed to be horizontal, we get $g_3=0$. Continuing for the forth coordinate in the basis $L_1,\ldots, L_n$, we obtain \begin{equation*}
\begin{split}
\dot\gamma(s) & =u_1L_1+u_2L_2+g_4L_4+(-2g_4c_1+g_5)\partial_5+\ldots+(-(n-3)g_4c_{n-4}+g_n)\partial_n.
\end{split}\end{equation*} To obtain the horizontal curve we take $g_4=0$.
Proceeding in the same way we conclude that a horizontal path satisfies the conditions (\ref{hc}).
\end{proof}
\s

\noindent{\it Remark.} Since we study left-invariant actions of $L_k$ on $\mathcal M_n$ we can take the vanishing initial conditions $c_k(0)=0$. So we may choose freely two coordinates $c_1$ and $c_2$ as two degrees of freedom. The resting
coordinates will be given as a solution to (\ref{hc}).
\s

\noindent{\it Remark.} Proposition \ref{hcprop} gives a  description of horizontal paths locally
in a neighborhood of the origin in $\mathcal M_n$. Checking the condition of horizontality (\ref{hc}) we must be sure
that the path lies inside $\mathcal M_n$. The L\"owner-Kufarev representation guarantees us this. For example,
for $n=3$, any L\"owner-Kufarev trajectory in $\mathcal M_3$ corresponding to an odd function $f(z)=z+c_2z^3+c_4z^5+\dots$ is horizontal. Just to make a concrete example, take the starlike function
\[
w(z,t)=\frac{e^{-t}z}{\sqrt{1+z^2(1-e^{-2t})}},
\]
with $p_1\equiv 0$, $p_2\equiv 1$, $p_3\equiv p_4\equiv\dots \equiv 0,$ and $c_1(t)\equiv 0$, $c_2(t)=\frac{1}{2}(e^{-2t}-1)$, $c_3(t)\equiv 0$, etc.

\subsection{Hamiltonian formalism for $\mathcal M_n$}

We choose the symplectic scalar product for $L_1,\dots, L_n$ to be given by the unit matrix $\{\delta_{j,k}\}$.
Being restricted onto the distribution $\mathcal D$ and taking into account the above matrix we get the Hamiltonian in the form 
$$
H(\xi_1,\ldots,\xi_n,c_1,\ldots,c_n)=|l_1|^2+|l_2|^2,$$ 
where
\begin{eqnarray*}
l_1&= &\bar\xi_1+2c_1\bar\xi_2+\ldots+n c_{n-1}\bar\xi_n,\\
l_2&=&\bar\xi_2+2c_1\bar\xi_3+\ldots+(n-1)c_{n-2}\bar\xi_n.
\end{eqnarray*}
Observe the similarity in formal variables were $\bar\psi_k=\partial_k=\bar\xi_k$ in (\ref{L1}, \ref{L2}).

The system of Hamiltonian equations is given by \begin{equation}
\label{hsys}\begin{split}
& \dot c_1  =\bar l_1\\
& \dot c_2  =2c_1\bar l_1+\bar l_2 \\
& \dot c_k=kc_{k-1}\bar l_1+(k-1)c_{k-2}\bar l_2,\quad k=3,\dots, n.\\
& \dot\xi_k =-(k+1)\xi_{k+1}l_1-(k+1)\xi_{k+2}l_2, \quad k=1,\dots, n-2, \\
& \dot\xi_{n-1}  = -n\xi_{n}l_1\\
& \dot\xi_{n} =0.\end{split}\end{equation}

\begin{proposition}\label{prop}
Any solution of the Hamiltonian system (\ref{hsys}) is a horizontal path.
\end{proposition}
\begin{proof}
Observe that 
\begin{equation}
\bar l_1=\dot c_1\quad\text{and}\quad \bar l_2= \dot c_2 -2c_1\dot c_1.\label{u=l}
\end{equation}
Substituting $\bar l_1$ and $\bar l_2$ into equations for $\dot c_3,\ldots,\dot c_n$, we obtain the horizontality conditions (\ref{hc}). 
 \end{proof}
 
 Likely for horizontal
paths, we assume vanishing initial conditions.

\begin{proposition}
Define $l_3$ as
\[
l_3=\bar \xi_3+2c_1\bar \xi_4+\ldots+(n-2)c_{n-3}\bar \xi_n.
\]
Then, 
\begin{itemize}

\item[(i)] $\dot l_1=\bar l_2 l_3$ and $\dot l_2=-\bar l_1 l_3$.

\item[(ii)] The energy of the system $\frac{1}{2}(|u_1|^2+|u_2|^2)$ is conserved along the geodesics.
The Carnot-Carath\'eodory length  of the tangent vector  is conserved along the geodesics.
\end{itemize}
\end{proposition}

\begin{proof}
The proof of (i) is straightforward. Differentiating $l_1$ and $l_2$ and using expressions for $\dot \xi_k$ and $\dot c_k$ from
the Hamiltonian system (\ref{hsys}) we obtain the necessary result. To prove (ii) we observe that
$\frac{\partial}{\partial t} (|l_1|^2+|l_2|^2)=0$ by (i). Moreover, the values of $u_1$ and $u_2$ coincide with $\bar l_1$ and $\bar l_2$ on geodesics by~(\ref{u=l}).
\end{proof}

As a consequence we get the solution to (\ref{hsys}) for $n=3$. Observe that $l_3=\bar\xi_3=const$ in this case. Hence,
$\ddot{c}_1=\dot{\bar{l}}_1=l_2\xi_3=\overline{(\dot c_2-2c_1\dot c_1)}\xi_3$ by the above proposition. We continue by $\ddot{c}_2=\frac{d^2}{d t^2}(c_1^2)+\dot{\bar l}_2=\frac{d^2}{d t^2}(c_1^2)-l_1\xi_3=\frac{d^2}{d t^2}(c_1^2)-\dot{ \bar c}_1 \xi_3$. Therefore, 
\begin{equation}
\label{eq2}
\begin{split}
&\ddot c_1+|\xi_3|^2c_1 =\bar K\xi_3,\\
&\dot c_2 =2c_1\dot c_1-\bar c_1\xi_3+K,
\end{split}
\end{equation} 
 where $K$ is a constant of integration and is calculated by
the initial speed $K=\dot c_2(0)$. The solution to the equation (\ref{eq2}) is
\[
c_1=Ae^{i|\xi_3| t}+B e^{-i|\xi_3| t}+{\bar K}/{\bar\xi_3}, \quad \mbox{where \ }A+B+{\bar K}/{\bar\xi_3}=0.
\]
Substituting $c_1$ in the equation for $c_2$ we get
\begin{eqnarray*}
c_2 &= &A^2e^{2i|\xi_3| t}+B^2e^{-2i|\xi_3| t}-2(Ae^{i|\xi_3| t}+Be^{-i|\xi_3| t})(A+B)\\& - &\frac{i\xi_3}{|\xi_3|}(Ae^{-i|\xi_3| t}-Be^{i|\xi_3| t}-(A-B))+4AB+A^2+B^2.
\end{eqnarray*}
The coordinate $c_3$ is calculated as a solution to the equation 
\[
\dot c_3= 3c_2\dot c_1+2c_1(\dot c_2-2c_1\dot c_1),\quad c_3(0)=0.
\]
The corresponding explicit expression is a matter of elementary calculations and we omit awkward formulas.

\s
\noindent{\it Remark.} Our Hamiltonian formalism and geodesics are linked to the sub-Riemannian geometry
led on $\mathcal M_n$ by Kirillov's vector fields. So there is no direct connection with the first
Hamiltonian system described in Section 2.3. The above Hamiltonian system (\ref{hsys}) gives local geodesics
in $\mathcal M_n$ about the origin and we do not expect any global description of geodesics because starting from the origin they may leave $\mathcal M_n$ in time.
\s

\subsection{Lagrangian formalism for $\mathcal M_3$}

Let us consider the Lagrangian function \begin{equation}\label{lagr}L(c,\bar c,\dot c,\bar{\dot c})= |\dot c_1|^2+|\dot c_2-2c_1\dot c_1|^2+\R\bar{\lambda}(\dot c_3-3c_2\dot c_1-2c_1\dot c_2+4c_1^2\dot c_1).\end{equation} It splits in two terms: the kinetic 
energy $|\dot c_1|^2+|\dot c_2-2c_1\dot c_1|^2$, and the  non-holonomic constraint $\dot c_3=3c_2\dot c_1+2c_1\dot c_2-4c_1^2\dot c_1$, that reflects the horizontality condition. We are interested in minimizing  the action integral $$S(c,\tau)=\int_{0}^{\tau}L(c,\bar c,\dot c,\bar{\dot c})\,ds.$$ The minimum of the action is attained at a critical curve $\zeta(s)$ satisfying the Euler-Lagrange system \begin{equation}\label{el}\frac{d}{ds}\Big(\frac{\partial L}{\partial\dot c}\Big)=\frac{\partial L}{\partial c},\quad \frac{d}{ds}\Big(\frac{\partial L}{\partial\dot{\bar c}}\Big)=\frac{\partial L}{\partial \bar c}.\end{equation}

\begin{proposition}\label{lem1}
The solution to the Euler-Lagrange system~\eqref{el} is a solution to the Hamiltonian system~\eqref{hsys} if and only if it is a horizontal path.
\end{proposition}

\begin{proof}
If the solution to the Euler-Lagrange system~\eqref{el} is a solution to the Hamiltonian system~\eqref{hsys}, then it is a horizontal path by Proposition \ref{prop}.

To show the reciprocal statement we perform  auxiliary calculation. Substituting the Lagrangian~\eqref{lagr} in the equations~\eqref{el}  we get \begin{equation}\label{sys1}
\begin{split}
\bar{\ddot c}_1-2c_1(\overline{\ddot c_2-\ddot{( c_1^2)}})-\bar\lambda\dot c_2 & =0\\
\overline{\ddot c_2-\ddot{( c_1^2)}}-\bar\lambda\dot c_1 & =0\\
\frac{d}{ds}\bar\lambda &=0.
\end{split}
\end{equation} We conclude that $\lambda$ is a constant.  Simplifying the first two equations we get 
\begin{equation}\label{sys2}
\begin{split}
\ddot c_1 & =\xi_3(\overline{\dot c_2-\dot{( c_1^2)}}) \\
\dot c_2 & =\dot{( c_1^2)}-\dot{\bar c}_1 \xi_3+K\\
\lambda & =\xi_3.
\end{split}
\end{equation} The latter equality is due to the Legendre transform. In the latter system we recognize the equations for geodesics (\ref{eq2}). 
\end{proof}

\subsection{Dual basis}

The following $1$-forms give the dual basis of the cotangent space for the basis $\{L_k\}$ of the tangent space:
\begin{equation}
\label{bf}\begin{split}
& \omega_1=dc_1\\
& \omega_2=dc_2-2c_1\omega_1\\
&\omega_k=dc_k-2c_1\omega_{k-1}-3c_2\omega_{k-2}-\ldots-kc_{k-1}\omega_1, \quad k=3,\dots, \infty. 
\end{split}\end{equation} We have $\omega_k(L_j)=\delta_{kj}$. Define the forms $\eta_k$ by 
\begin{equation}
\label{cf}
\eta_k=dc_k-kc_{k-1}\omega_{1}-(k-1)c_{k-2}\omega_{2}\quad k=3,\dots, \infty.
\end{equation} 
Then the form $\eta=\sum_{k=1}^n$ defines the distribution $\mathcal D$ for $\mathcal M_n$ as a kernel $$\mathcal D=\{X\in T\mathcal M_n:\ \eta(X)=0\}.$$  

A {\it contact form} $\alpha$ on a $(2n+1)$-dimensional manifold is a local $1$-form with the property $$\alpha\wedge d \alpha\neq 0.$$

 In our case $n$ is the complex dimension. Nevertheless, for $n=3$ we get
\[  
\eta_3\wedge d\eta_3 =dc_1\wedge dc_2\wedge dc_3\neq 0.
\]
The form $\eta_3$ is contact and its kernel defines the distribution $\mathcal D$ in $\mathcal M_3$.


\begin{thebibliography}{99}

\bibitem{Airault} H.~Airault, P.~Malliavin, {\it Unitarizing probability measures for representations of Virasoro
algebra}, J. Math. Pures Appl. {\bf 80}
 (2001), no. 6, 627--667.
 
\bibitem{Arnold} V.~I.~Arnold, {\it Mathematical methods of classical
mechanics}, Springer-Verlag, New York, 1989

\bibitem{Babelon} O.~Babelon, D.~Bernard, M.~Talon, {\it Introduction to
classical integrable systems}. Cambridge Monographs on Mathematical
Physics. Cambridge University Press, Cambridge, 2003.

\bibitem{Babenko}
K.~I.~Babenko, {\it The theory of extremal problems for univalent
functions of class $S$}, Proc. Steklov Inst. Math., No. 101 (1972).
Transl. American Mathematical Society, Providence, R.I., 1975.

\bibitem{Bambusi} D.~Bambusi, G.~Gaeta, {\it  On persistence of invariant tori and a theorem by Nekhoroshev},
Math. Phys. Electron. J. 8 (2002), Paper 1, 13 pp.

\bibitem{Bolsinov} A.~V.~Bolsinov, A.~T.~Fomenko, {\it Integrable Hamiltonian systems.
Geometry, topology, classification}. Chapman \& Hall/CRC, Boca
Raton, FL, 2004.

\bibitem{Bott}
R.~Bott, {\it On the characteristics classes of groups of diffeomorphisms}, Enseignment Math. (2) {\bf 23} (1977), 209--220.


\bibitem{Branges} L.~de~Branges, {\it A proof of the Bieberbach conjecture}, Acta
Math. {\bf 154} (1985), no. 1-2, 137--152.

\bibitem{Chow}
W.~L.~Chow, {\it Uber Systeme Von Lineaaren Partiellen Differentialgleichungen erster Ordnung}, Math. Ann. {\bf 117} (1939), 98--105.

\bibitem{Duren}
P.~Duren, {\it Univalent functions}, Springer, New York, 1983.

\bibitem{Fiorani} E.~Fiorani, G.~Giachetta, G.~Sardanashvily, {\it The Liouville-Arnold-Nekhoroshev
theorem for non-compact invariant manifolds}, J. Phys. A {\bf 36}
(2003), no. 7, L101--L107.

\bibitem{Gaeta} G.~Gaeta, {\it The Poincar\'e-Lyapounov-Nekhoroshev
theorem}, Ann. Physics {\bf 297} (2002), no. 1, 157--173.

\bibitem{Gelfand}
I.~M.~Gel'fand, D.~B.~Fuchs, {\it Cohomology of the Lie algebra of vector
fields on the circle},   Functional Anal. Appl. {\bf 2} (1968), no.4, 342--343.

\bibitem{Hoermander}
L.~H\"ormander, {\it Hypoelliptic second order differential equations},  Acta Math.  {\bf 119}  (1967), 147--171.

\bibitem{Kir}
A.~A.~Kirillov, {\it K\"ahler structure on the $K$-orbits of a
group of diffeomorphisms of the circle}, Functional Anal. Appl. {\bf 21} (1987), no. 2, 122--125.

\bibitem{Kir2}
A.~A.~Kirillov, {\it Geometric approach to discrete series of unirreps for Vir}, J. Math. Pures Appl. {\bf 77} (1998),  735--746.


\bibitem{KY1}
A.~A.~Kirillov, D.~V.~Yuriev, {\it K\"ahler geometry of the
infinite-dimensional homogeneous space $M={\rm Diff}\sb +(S\sp
1)/{\rm Rot}(S\sp 1)$}, Functional Anal. Appl. {\bf 21} (1987), no. 4, 284--294.

\bibitem{KY2}
A.~A.~Kirillov, D.~V.~Yuriev, {\it Representations of the Virasoro
algebra by the orbit method}, J. Geom. Phys. {\bf 5} (1988), no.
3, 351--363.

\bibitem{Mitchell}
J.~Mitchell, {\it On Carnot-Carath\'eodory metrics},  J. Differential Geom. {\bf 21}  (1985),  no. 1, 35--45. 

\bibitem{Montgomery}
R.~Montgomery, {\it A survey of singular curves in sub-Riemannian geometry}, J. Dynamical and Contr. Syst. {\bf 1} (1995), no. 1, 49--90.

\bibitem{NSW}
Nagel~A., Stein~E.~M., Wainger~S. {\it Balls and metrics defined
by vector fields. I. Basic properties.}  Acta Math. \textbf{155}
(1985), no. 1-2, 103--147.


\bibitem{Nekhoroshev}  N.~N.~Nekhoroshev, {\it  The Poincar\'e-Lyapunov-Liouville-Arnol'd
theorem},  Functional Anal. Appl.  {\bf 28}  (1994),
no. 2, 128--129.

\bibitem{Neretin} Yu.~A.~Neretin, {\it Representations of Virasoro and affine Lie algebras},
Encyclopedia of Mathematical Sciences, Vol. 22, Springer-Verlag, 1994, pp. 157--225.

\bibitem{Pansu}
P.~Pansu, {\it M\'etriques de Carnot-Carath\'eodory et quasiisom\'etries des espaces sym\'etriques de rang un},  Ann. of Math. (2)  {\bf 129}  (1989),  no. 1, 1--60.

\bibitem{Pom}
Ch.~Pommerenke, {\it Univalent functions, with a chapter on
quadratic differentials by G.~Jensen}, Vandenhoeck \& Ruprecht,
G\"ottingen, 1975.

\bibitem{Prokhorov}
D.~Prokhorov,  {\it Sets of values of systems of functionals in
classes of univalent functions}, Mat. Sb. {\bf 181} (1990), no. 12,
1659--1677; translation in Math. USSR-Sb. {\bf 71} (1992), no. 2,
499--516.

\bibitem{ProkhVas}
D.~Prokhorov, A.~Vasil'ev, {\it Univalent functions and integrable systems}, Commun. Math. Phys. {\bf 262} (2006), no. 2,
393--410.

\bibitem{Rashevski}
P.~K.~Rashevski, {\it About connecting two points of complete nonholonomic space by admissible curve}, Uchen. Zap. Ped. Inst. K.~Libknehta (1938), 83--94.

\bibitem{SS}
A.~C.~Schaeffer, D.~C.~Spencer, {\it Coefficient Regions for
Schlicht Functions (With a Chapter on the Region of the Derivative
of a Schlicht Function by Arthur Grad)}, American Mathematical
Society Colloquium Publications, Vol. 35. American Mathematical
Society, New York, 1950.

\bibitem{Strichartz1}
R.~S.~Strichartz, {\it Sub-Riemannian geometry}, J.~Differential Geom. {\bf 24} (1986), no. 2, 221--263.

\bibitem{Strichartz2}
R.~S.~Strichartz, {\it Corrections to: ``Sub-Riemannian geometry"}, J.~Differential Geom. {\bf 30} (1989), no. 2, 595--596.

\bibitem{Vas}
A.~Vasil'ev, {\it Energy characteristics of subordination chains}, arXiv: math-ph/0509072, 2005, 12 pp.

\bibitem{Zakharov} ed V.~E.~Zakharov, {\it What is integrability?}, Springer Series
in Nonlinear Dynamics. Springer-Verlag, Berlin, 1991.

\end{thebibliography}
\end{document}